\documentclass[11pt,leqno]{article}
\date{}
\usepackage{amssymb,amsbsy,amsmath,amsfonts,amssymb,amscd, mathrsfs}
\usepackage[english]{babel}
\usepackage[T1]{fontenc}
\usepackage{indentfirst}
\usepackage{color}
\usepackage{authblk}

\makeatletter
\@addtoreset{equation}{section}
\makeatother

\newtheorem{statement}{}[section]
\newtheorem{theoreme}[statement]{Theorem}
\newtheorem{lemme}[statement]{Lemma}

\newtheorem{proposition}[statement]{Proposition}

\newcommand\C{\mathbb C}

\newcommand\R{\mathbb R}
\newcommand\T{\mathbb T}
\newcommand\D{\mathbb D}
\newcommand\Z{\mathbb Z}

\newcommand\ind{{\rm 1\kern-.30em I}}

\renewcommand \Re{{\mathfrak R}{\rm e}\,}

\let\phi=\varphi
\newenvironment{demo}{\medbreak\begingroup\noindent{\bf Proof : }}{\hfill$\square$\endgroup\goodbreak\medskip}

\title{\bf On Bernstein's inequality for polynomials}

\author[1]{H.~Queff\'elec\thanks{Herve.Queffelec@univ-lille.fr}}
\author[2]{R. Zarouf\thanks{rachid.zarouf@univ-amu.fr}}
\affil[1]{Universit\'{e} Lille Nord de France, USTL, Laboratoire Paul Painlev\'{e} U.M.R. CNRS 8524 et F\'{e}d\'{e}ration CNRS Nord-Pas-de-Calais FR 2956 F-59 655 Villeneuve d'Ascq Cedex, France.}
\affil[2]{Aix-Marseille Universit\'{e}, Laboratoire Apprentissage, Didactique, Evaluation,
Formation, 32 Rue Eug\`{e}ne Cas CS 90279 13248 Marseille Cedex 04, France}
\affil[2]{Department of Mathematics and Mechanics, Saint Petersburg State University,
28, Universitetski pr., St. Petersburg, 198504, Russia}

\date{\footnotesize \today}

\begin{document}

\maketitle
\begin{abstract}

Bernstein's classical inequality asserts that given a trigonometric polynomial $T$ of degree $n\geq1$, the sup-norm of the derivative of $T$ does not exceed $n$ times the sup-norm of $T$.  We present various approaches to prove this inequality and some of its natural extensions/variants, especially when it comes to replacing the sup-norm with the $L^p-norm$.

\end{abstract}

\vspace{5mm}
\section{Introduction}
\noindent Bernstein's inequality for trigonometric polynomials (\cite{Bernstein1}), already one century old,  played a fundamental role in harmonic and complex Analysis, as well as in approximation theory (\cite{Bernstein1}, \cite{Bernstein2}, \cite{LOR}) and in the study of random trigonometric series (\cite{SZ}, \cite{KAHA} Chapter 6) or random Dirichlet series (\cite{QU}, Chapter 5), when generalized to several variables in the latter case. One can also mention its use in the theory of Banach spaces (\cite{PI}, p. 20-21), or its extensive use in Numerical Analysis.\\

\noindent The purpose of this survey is not to focus on applications of Bernstein's inequality, but on various approaches (some classical, some more recent) to proving this inequality and its extensions, and to show that even if it was  first stated for the sup-norm, it  is valid for a {\bf large class of norms, even of quasi-norms}. We will not be interested in describing all equality cases. There will be two key words here:\\

\noindent $\bullet$ Convexity (with both  real and complex variable approaches), well adapted to norms.\\

\noindent $\bullet$ Subharmonicity (with a rather complex variable approach), better adapted to quasi-norms.
\bigskip

The paper is organized as follows:\\

\noindent  $\bullet$ Section 1 is this introduction, with reminders and the proof of Riesz.\\

\noindent  $\bullet$ Section 2, essentially a real variable section, illustrates the role of convexity and of translation-invariance in generalized forms of Bernstein's inequality for Fourier transforms of compactly supported measures.\\

\noindent  $\bullet$ Section 3 is a transition between real and complex methods.\\

\noindent  $\bullet$ Section 4 is a section using intensively complex and hilbertian methods (integral representations, reproducing kernels), and new Banach algebra norms (Wiener norm,  Besov norm) in connection with operator theory and functional calculus. Embedding inequalities other than Bernstein's one will also be considered.\\

\noindent  $\bullet$ Section 5\  ``jumps'' into quasi-norms with a somewhat extreme case, the Mahler $\Vert .\Vert_0$ quasi-norm, called Mahler norm for simplicity.\ Here, subharmonicity plays a key role.\\

\noindent  $\bullet$ Section 6\    ``climbs again the road'' from the quasi-norm $\Vert .\Vert_0$ to quasi-norms or norms $\Vert .\Vert_p$ with $0<p\leq \infty$, through integral representations. \\

\noindent  $\bullet$ The final Section 7 concludes with some remarks and open questions.

\subsection{Reminders and notations}
Bernstein's inequality is generally quoted under the following form.
\begin{theoreme}\label{sb} Let $T(x)=\sum_{k=-n}^n a_k e^{ikx}$ be a trigonometric polynomial of degree $\leq n$. Then:
$$\sup_{x\in \R}\vert T' (x)\vert \leq n\sup_{x\in \R}\vert T(x)\vert$$
and the constant $n$ is optimal in general ($T(x)=e^{inx}$).
\end{theoreme}
Throughout this paper, we shall have to make a careful distinction between  trigonometric polynomials as above and algebraic polynomials
$$T(x)=\sum_{k=0}^n a_k e^{ikx}=P(e^{ix}) \hbox{\ with}\ T'(x)=ie^{ix}P'(e^{ix}) \hbox{\ and}\ |T'(x)|=|P'(e^{ix})| $$ where $P$ is the ordinary polynomial $P(z)=\sum_{k=0}^n a_k z^k$ for which complex techniques are more easily available. If  once and for all   $\D$ designates the open unit disk  and $\T=\{z: |z|=1\}$ its boundary, as well as $\Vert f\Vert_\infty=\sup_{z\in \D} |f(z)|$ when $f$ is a bounded analytic function on $\D$, the maximum modulus principle gives us for $T(x)=P(e^{ix})$ as above:
$$\Vert T\Vert_\infty=\sup_{x\in \R}|P(e^{ix})|=\sup_{z\in \D}|P(z)|=\Vert P\Vert_\infty,$$
and we shall always identify both sup-norms, as well as $P$ and $x\mapsto P(e^{ix})$. \\
The Haar measure of $\T$ will be denoted $m$:
$$\int_{\T}fdm=\int_{0}^1 f(e^{2i\pi\theta})d\theta.$$
The $L^p$-norm (quasi-norm when $0<p<1$) will always refer to the measure $m$. We also set (the Mahler norm)
\begin{equation}\label{ma} \Vert f\Vert_0=\lim_{p\to 0} \Vert f\Vert_p=\exp\big(\int_{\T} \log |f|dm\big).\end{equation}
Recall that  $\Vert f\Vert_\infty=\lim_{p\to \infty} \Vert f\Vert_p$.\\
\medskip

\noindent If $P(z)=\sum_{k=0}^n a_k z^k=a\prod_{j=1}^n (z-z_j)$ is an algebraic polynomial, its (complex) reciprocal polynomial $Q$ is defined by
\begin{equation}\label{re} Q(z)=z^n \overline{P(1/\overline{z})}=\sum_{k=0}^n \overline{a_{n-k}}\, z^k=\overline{a}\prod_{j=1}^n (1-\overline{z_j}z).\end{equation}
The reciprocal polynomial of $Q$ is $P$.The following obvious property is quite useful:
\begin{equation}\label{ob} |z|=1\Rightarrow  |Q(z)|=|P(z)|.\end{equation}
\medskip
 For $P$ as above, Jensen's formula tells that
\begin{equation}\label{jf}\Vert P\Vert_0=|a|\prod_{j=1}^n \max(1,|z_j|).\end{equation}
\subsection{Bernstein through interpolation, \ Riesz formula}
Bernstein  (\cite{Bernstein1}) initially obtained $$\Vert T' \Vert_\infty \leq 2n\Vert T\Vert_\infty$$ and the best constant $n$ was shortly afterwards obtained by  E.~Landau (\cite{Bernstein3}) by a reduction to a sum of sines, and slightly later by  M.~Riesz (\cite{Riesz}), using a new interpolation formula.\\
  We first present the proof of Riesz. See also the nice books \cite{PRA} page 146 and \cite{BOER} page 178.
\begin{theoreme}[M.Riesz.]There exist $c_1,\ldots, c_{2n}\in \C$ and $x_1,\ldots, x_{2n}\in \R$ with $\sum_{r=1}^{2n} |c_r|=n$ such that, for all trigonometric polynomials $T$ of degree $n$:
\begin{equation}\label{intfin}T'(x)=\sum_{r=1}^{2n} c_r T(x+x_r) \hbox{\quad for all}\  x\in \R. \end{equation}
In particular
$$|T'(0)| \leq n \sup_{x\in E_n}|T(x)|$$
where $E_n=\{x_j, 1\leq j\leq 2n\}$.
\end{theoreme}

\begin{demo} we  sketch the proof. Let  $r$ be an integer with  $1\leq r\leq 2n$. We set:
\begin{equation}\label{nota}x_r=\frac{(2r-1)\pi}{2n},\quad  \omega=e^{\frac{i\pi}{2n}},\quad  z_r=e^{ix_r}=\omega^{2r-1},\quad z_{r}^{2n}=-1. \end{equation}
\noindent  An easy variant of the Lagrange interpolation formula for the   $2n$ points $z_r$ gives, for any polynomial  $P(z)=\sum_{k=0}^{2n}c_k z^k$:
\begin{equation}\label{lagrangest} P(z)=\frac{z^{2n}+1}{2}\,(c_0+c_{2n})+\frac{z^{2n}+1}{2}\frac{1}{2n}\sum_{r=1}^{2n} P(z_r)\frac{z_r+z}{z_r -z} \cdot\end{equation}
Next, if  $T(x)=\sum_{k=0}^n \big(a_k \cos kx+b_k \sin kx\big)$,  we apply  formula (\ref{lagrangest}) to the polynomial  $P(z)=\sum_{k=0}^{2n}c_k z^k$  defined by
$P(e^{ix})=e^{inx}\, T(x).$  We get
\begin{equation}\label{marcelo}T(x)=a_n \cos nx+\cos nx\,\frac{1}{2n} \sum_{r=1}^{2n} T(x_r)(-1)^{r+1}\cot \frac{(x_r- x)}{2} \cdot\end{equation}
  Differentiation at $0$ now gives
\begin{equation}\label{marcelor}T'(0)=\frac{1}{2n} \sum_{r=1}^{2n} T(x_r)\frac{(-1)^{r+1}}{2\sin^{2}(x_r/2)}=:\sum_{r=1}^{2n} c_r T(x_r)\hbox{\ with}\ \sum_{r=1}^{2n} |c_r|=n.\end{equation}
By translation, we get \textit{the Riesz interpolation formula:}
\begin{equation}\label{intfin}T'(x)=\sum_{r=1}^{2n} c_r T(x+x_r).\end{equation}
In convolution terms:

\begin{equation} \label{mr}T'=T\ast \mu_n, \hbox{\quad where}\quad \mu_n=\sum_{r=1}^{2n} c_r \delta_{x_r} \hbox{\ and}\ \Vert \mu_n\Vert=n\end{equation}
 and this clearly ends  the  proof.\end{demo}
 \noindent  Observe that the measure $\mu_n$ is a \textit{\ finite} combination of Dirac point masses. We will later see an extension of this method in which $\mu_n$ is discrete, but an infinite combination of Dirac point masses.

\section{Convexity}
We begin with giving a general form (due to R.~P. Boas) of Bernstein's previous inequality, as in the book \cite{KAH} page 30, and which may be seen as an extension of Riesz's proof. This form is valid for non-periodic (almost periodic) trigonometric polynomials as well. We denote the derivative $f'$ by $Df$ and the translate of $f$ by a real number $a$ by $T_{a}f$, that is $T_{a}f(x)=f(x+a)$. The convolution of the function $f$ and the measure $\mu$ (already appearing in Riesz's proof) is accordingly defined as
$$f\ast \mu=\int_{\R} (T_{t}f) d\mu(t).$$
\bigskip

\begin{theoreme}\label{rp} Let $\lambda>0$.Then there exists a complex sequence $(c_k)_{k\in \Z}$ and a real sequence   $(t_k)_{k\in \Z}$, depending only on $\lambda$,  such that $\sum_{k\in \Z}|c_k|=\lambda$ and that, whenever
 $f(x)=\int_{\R} e^{itx}d\mu(t)$ is the Fourier transform of a complex measure on $\R$ supported by $[-\lambda,\lambda]$, then
$$Df=\sum_{k\in \Z} c_k T_{t_k}f.$$
If one prefers, $Df=f\ast \mu$ with $\Vert \mu\Vert=\lambda$.\\
In particular, if $f(x)=\sum_{j=1}^N a_j e^{i\lambda_{j}x}$ where the $\lambda_j$'s are real and distinct with $|\lambda_j| \leq \lambda$, then
$$\Vert f' \Vert_\infty \leq \lambda\Vert f\Vert_\infty.$$

\end{theoreme}
\begin{demo}we rely on the following lemma.
 \begin{lemme}\label{zg} Let $\varphi$ be the $4\lambda$-periodic odd function defined by

\[
\varphi(t)=\begin{cases}
\:t & \mbox{\rm{if}}\ 0\leq t\leq\lambda\\
\:2\lambda-t & \mbox{\rm{if}}\ \lambda\leq t\leq2\lambda .\\
\end{cases}
\]

Then,
$$i\varphi(t)=\sum_{k\in \Z} c_k e^{i\pi kt/ 2\lambda} \hbox{\ with}\ \sum_{k\in \Z} |c_k|=\lambda.$$
\end{lemme}
Indeed, let us work with the space $E$ of $4\lambda$-periodic functions, initially defined on $[-2\lambda,2\lambda]$. Let $\chi\in E$ be the characteristic function of the interval $[-\lambda,\lambda]$. Let $\psi(t)=\varphi(t+\lambda)+\lambda\in E$, a triangle function on $[-2\lambda, 2\lambda]$.  We see that $\psi=4\lambda(\chi\ast\chi)$,  and it can hence be written as
$\psi(t)=\sum_{k\in \Z} d_k e^{i\pi kt/ 2\lambda}$ with $d_k=4\lambda(\widehat{\chi}(k))^2\geq 0$, so that  $\sum_{k\in \Z} d_k=\psi(0)=2\lambda$ and $d_0=\frac{1}{4\lambda}\int_{-2\lambda}^{2\lambda} \psi(t)dt=\lambda$. Since $\varphi(t)=\psi(t-\lambda)-\lambda$, the lemma follows with $c_0=0$ and $c_k=id_k i^{-k}$ if $k\neq 0$.
\medskip

Coming back to Theorem \ref{rp}, we see that, since $\mu$ is supported by $[-\lambda, \lambda]$:
$$f'(x)=\int_{\R} it e^{itx}d\mu(t)=\int_{\R} i\varphi(t) e^{itx}d\mu(t)=\sum_{k\in \Z} c_k \int_{\R}e^{i\pi kt/2\lambda}e^{itx}d\mu(t)$$$$=\sum_{k\in \Z} c_k f(x+t_k) \hbox{\ with}\  t_k=\frac{k\pi}{2\lambda}$$
and this ends the proof of the general part of our theorem. For the special case, just observe that $f$ is the Fourier transform of the discrete measure $\mu=\sum_{j=1}^N a_j \delta_{\lambda_j}$. The meaning of this theorem is that, for Fourier transforms of compactly supported measures, the differential operator $D$ can be replaced by kind of a convex combination of translation operators $T_{t_k}$; this is why Theorem \ref{rp} belongs to convexity.\end{demo}
\medskip

\noindent It is worth mentioning a more general application of Theorem \ref{rp}.
\begin{theoreme}\label{pw}\label{co} Let $f$ be an entire function of exponential type $\lambda$ (namely $|f(z)|\leq Ce^{\lambda |z|}$), bounded on the real axis ($\Vert f\Vert_\infty:=\sup_{x\in\R} |f(x)|<\infty$).  Then
$$\Vert f' \Vert_\infty \leq \lambda\Vert f\Vert_\infty.$$
\end{theoreme}
\begin{demo}indeed, by the Paley-Wiener theorem (\cite{KAT} page 212), $f$ restricted to the real line is the Fourier transform of a measure (indeed of an $L^2$-function) supported by $[-\lambda,\lambda]$.\end{demo}

\noindent Let us denote by $\mathcal{P}_n$ the translation-invariant space of trigonometric polynomials $\sum_{|j|\leq n} p_j e^{ijx}$ of degree $\leq n$. A nice corollary of Theorem \ref{rp} is the following:

\begin{theoreme}\label{co} Let $\Vert.\Vert$ be a \textnormal{translation-invariant} norm on  $\mathcal{P}_n$. Then
$$\Vert f' \Vert \leq n\Vert f\Vert \hbox{\ for all} \ f\in \mathcal{P}_n.$$
\end{theoreme}
\begin{demo}writing $f'=\sum_{k\in \Z} c_k T_{t_k}f$ and taking norms (note that the series on the right-hand side is absolutely convergent for the norm  $\Vert.\Vert$), we get
$$\Vert f'\Vert\leq \sum_{k\in \Z} |c_k|\,\Vert T_{t_k}f\Vert = \sum_{k\in \Z} |c_k|\,\Vert f\Vert =n\Vert f\Vert.$$
\end{demo}
This can be applied to the $L^p$-norm with respect to the Haar measure $m$ of the circle $\T$, with $1\leq p\leq \infty$, more generally to the $L^\psi$-norm where $\psi$ is any Orlicz function (\cite{ZYG} page 173).  There are lots of applications, and improvements of the factor $n$ under special assumptions (as unimodularity of coefficients);  we just mention the paper \cite{QUSA} and the book \cite{KAHA} with applications to random Fourier series.\\

\noindent As we will now see, subharmonicity and complex methods allow us to go beyond convexity and to consider $L^p$-quasi-norms for $0<p<1$, even for $p=0$ (the Mahler norm). We begin with a ``transition'' section.
\medskip
\section{Convexity and Complexity}
What follows still belongs to convexity, in spite of the appearance of Complex Analysis and the maximum principle, behind which subharmonicity is lurking. Let us consider this section as a transition, we will be more explicit on subharmonicity later. A typical example of this transition is the famous Gauss-Lucas theorem, and its extension by Laguerre.

 \begin{theoreme}\label{gl} Let $f$ be an algebraic polynomial of degree $n$, all of which roots lie in a convex set $K$ of the plane. Then, all  the roots of the derivative $f'$ also lie in $K$.
 \end{theoreme}
 The following variant, due to Laguerre, of Theorem \ref{gl} is worth mentioning, in view of the forthcoming applications.
 \begin{theoreme}\label{la} Let $\rho\geq 1$.  Let $P$ be an algebraic polynomial of degree $n$, all of which roots lie inside $E:=\{z: |z|\geq \rho\}$.  Assume that $\xi, z$ satisfy
 $$(\xi-z)P'(z)+nP(z)=0.$$
 Then,  either $\xi$ or $z$   lie in $E$. \\
 As a consequence, $$|z|=1\Rightarrow \rho |P'(z)|\leq |Q'(z)|$$
 where $Q$ is the (complex) reciprocal polynomial of $P$.
 \end{theoreme}
\begin{demo} without loss of generality, we can assume that $\rho>1$.
 Denote here by $T_z$ the ``inversion'' with pole $z$, namely
 $$T_{z}(u)=\frac{1}{u-z}\cdot$$
 Let $F=\C\setminus E$ and  $z_1,\ldots, z_n$ the roots of $P$. In view of the formula $P'(z)/P(z)=\sum_{j=1}^n 1/(z-z_j)$, the relation between $z$ and $\xi$ can be written as
 \begin{equation}\label{bar}T_{z}(\xi)=\frac{1}{n}\sum_{j=1}^n T_{z}(z_j).\end{equation}
 We see that, modulo $T_z$, $\xi$ is none other than the barycenter of the $z_j$'s and convexity is again implied. Suppose now that $z\in F$. Then, $\infty=T_{z}(z)\in T_{z}(F)$, hence $T_{z}(F)$ is unbounded, and its complement $T_{z}(E)$ is {\bf convex} since it is a disk or a half-plane. Now, (\ref{bar}) shows that $T_{z}(\xi)\in T_{z}(E)$ since $z_j\in E, 1\leq j\leq n$, by hypothesis. That is $\xi\in E$. Finally, fix $z$ with modulus one.  Note that, since then $$z\frac{P'(z)}{P(z)}=n-\frac{\overline{zQ'(z)}}{\overline{Q(z)}},$$
 we have as well
 \begin{equation}\label{asw}\xi=-n\frac{P(z)}{P'(z)}+z= -\frac{\overline{zQ'(z)}}{\overline{Q(z)}} \times \frac{P(z)}{P'(z)}\cdot \end{equation}
 Now, $z\notin E$ since $\rho>1$. The first part of the theorem gives $\xi\in E$ or again $|\xi|\geq \rho$, giving the conclusion in view of (\ref{asw}) and of $|P(z)|=|Q(z)|$. \end{demo}
\medskip

\noindent The key point  of the end of this section is the following lemma of  term by term differentiation of inequalities; here, two  polynomials are involved:
\bigskip

 \begin{lemme}\label{duo} Let $f, F$ be two algebraic polynomials of degree $n$ satisfying \begin{enumerate}

\item $|f(z)|\leq |F(z)|$ for all $z\in \T$;
\item all roots of $F$ lie in the closed disk $\overline{\D}$.

\end{enumerate}
Then \begin{enumerate}
\item $|f(z)|\leq |F(z)|$ for all $|z|\geq 1$;
\item $|f'(z)|\leq |F'(z)|$ for all $z\in \T$.

\end{enumerate}
\end{lemme}

\begin{demo} We begin with  assertion $1$. Suppose first that  all roots of  $F$ lie in $\D$.  We consider the rational function $f/F$ in the (unbounded) open set  $\Omega=\{z : |z|>1\}$; this function is holomorphic and bounded in $\Omega$ (since $f$ and $F$ have the same degree  $n$) and  continuous on  $\overline{\Omega}$ since by hypothesis all roots of $F$ lie in $\D$. Moreover, $f/F$ has modulus $\leq 1$ on $\partial{\Omega}$. The  maximum modulus principle gives the conclusion. In the general case, one writes  (note that the multiplicity of the zero $z_j$ is higher for $f$ than for $F$, due to our first assumption)
$$ f(z)=\prod_{|z_j|=1}(z-z_j)^{\alpha_j} g(z) \hbox{\ and}\  F(z)=\prod_{|z_j|=1}(z-z_j)^{\alpha_j} G(z)$$
where $g$ and $G$ are  polynomials of the same degree,  all roots of $G$ lying in  $\D$, and satisfying  $|g(z)|\leq |G(z)|$ for $z\in \T$. From the first case, one gets
$$|z|\geq 1\Rightarrow |g(z)|\leq |G(z)|\Rightarrow |f(z)|\leq |F(z)|.$$

\noindent  \textit{The second assertion follows from the first one.} Let us indeed fix a complex number $w$ with modulus $>1$. If $|z|> 1$, we have
$$|w F(z)-f(z)|\geq |w|\,|F(z)|-|f(z)|\geq (|w|-1)|F(z)|>0,$$ and all the roots of the polynomial $w F-f$ lie in $\overline{\D}$, as well as  (by Theorem \ref{gl}) those of the derivative  $w F'-f'$. In particular:
$$|z|>1\Rightarrow f'(z)\neq w F'(z).$$
 By the  Gauss-Lucas theorem again, we have $F'(z)\neq 0$, therefore $f'(z)/F'(z)\neq w$.
The quotient $f'(z)/F'(z)$ being different from any   complex number of  modulus $>1$, we get:
$$|z|>1\Rightarrow |f'(z)|\leq  |F'(z)|.$$
Letting $|z|$ tend to $1$ gives the claimed result. \end{demo}
\medskip

\noindent {\bf Remark.} The previous lemma contains Bernstein's inequality for algebraic polynomials (meaning $f(e^{it})=\sum_{k=0}^n a_k e^{ikt}$), assuming that $|f(z)|\leq 1$ for $z\in \T$ and taking then $F(z)=z^n$. But the extension to trigonometric polynomials is not straightforward, and will need the full generality of Lemma \ref{duo}, under the following form (\cite{MAL}).

\begin{theoreme}[Malik.]\label{ma} Let $P$ be an algebraic polynomial of degree  $n$, and  $Q$ its  reciprocal polynomial.
We assume that  $\Vert P\Vert_\infty\leq 1$. Then
$$z\in \T\Rightarrow |P'(z)|+|Q'(z)|\leq n.$$

\end{theoreme}
\bigskip
\begin{demo} let $|w|>1$. It suffices to apply Lemma \ref{duo}  to $f=Q-\overline{w}$ and its reciprocal polynomial $F=P-w z^n$, which satisfy:   $|f|=|F|$ on $\T$, and   $F$ has no zeros outside $\D$,  by a new application of Lemma \ref{duo} to $P$ and $z^n$. We get for $z\in \T$:
$$|Q'(z)|\leq |P'(z)-w nz^{n-1}|,$$
whence the result by  adjusting the argument of $w$ and by letting its modulus tend to $1$. \end{demo}
\smallskip

\noindent Lax proved that if an algebraic polynomial $P$ of degree $n$ has no roots in $\D$, Bernstein's inequality can be improved as follows:   $\Vert P'\Vert_\infty\leq \frac{n}{2}\Vert P\Vert_\infty,$ the inequality being optimal. What precedes provides a simple proof and extension of Lax's result, due to Malik (\cite{MAL}).
\begin{theoreme} Let $\rho\geq 1$ and let $P$ be an algebraic polynomial f degree $n$, all of which roots have modulus $\geq \rho$. Then
$$\Vert P'\Vert_\infty \leq \frac{n}{1+\rho} \Vert P\Vert_\infty.$$
The constant $\frac{n}{1+\rho}$ is optimal.
\end{theoreme}
\begin{demo}the optimality is clear by considering $P(z)=\Big(\frac{z+\rho}{1+\rho}\Big)^n\cdot$ Now, assume that $\Vert P\Vert_\infty=1$ and fix a unimodular complex number $z$.  We combine two previous results:
$$ |P'(z)|+|Q'(z)|\leq n $$
$$ \rho|P'(z)|\leq |Q'(z)|\leq n.$$
We hence get
$$(1+\rho)|P'(z)|\leq n.$$
This ends the proof.\end{demo}
It is worth mentioning a corollary of Lax-Malik's result, due to Ankeny and Rivlin for $\rho=1$.
\begin{proposition} Let $\rho\geq 1$ and  $P$ be an algebraic polynomial of degree  $n$ with no roots in $\rho\D$. Then:
$$|z|\geq 1\Rightarrow |P(z)|\leq \frac{|z|^n +\rho}{1+\rho}\,\Vert P\Vert_\infty.$$
\end{proposition}
\begin{demo}we can assume $\Vert P\Vert_\infty=1$. By Malik's theorem, one gets that $|P'(z)|\leq  \frac{n}{1+\rho}$ for $|z|=1$. By  the maximum principle,
$|P'(z)|\leq  \frac{n}{1+\rho} |z|^{n-1}$ for $|z|\geq1$. Now, if $R>1$ and $\theta\in \R$, one can write
$$P(Re^{i\theta})-P(e^{i\theta})=\int_{1}^R e^{i\theta}P'(re^{i\theta})dr$$
whence
$$\big|P(Re^{i\theta})-P(e^{i\theta})\big|\leq \int_{1}^R \frac{n}{1+\rho}r^{n-1}dr=\frac{R^n -1}{1+\rho}\cdot$$
The triangle inequality now gives
$$\big|P(Re^{i\theta})\big|\leq\frac{R^n +\rho}{1+\rho}\cdot$$
This ends the proof.
\end{demo}

\noindent  Here is an interesting variant, and strenghtening, of Bernstein's inequality.

\begin{theoreme}[Schaake-van der Corput.] Let $T$ be a  \textnormal{real}   trigonometric polynomial of  degree $n$, with $|T(x)|\leq 1$ for all $x\in \R$. Then
$$(T'(x))^2+n^2 (T(x))^2\leq n^2 \hbox{\ for all}\ \theta \in \R.$$
In particular, $|T'(x)|\leq n$.

\end{theoreme}
 \bigskip
\begin{demo}let $P(e^{ix})=e^{inx} T(x)$, an   algebraic polynomial of degree  $2n$, and  $Q$ be its reciprocal polynomial.  Since   $T$ is real, we have:
$$Q(e^{ix})=e^{2inx}\overline{P(e^{ix})}=e^{2inx}e^{-inx}T(x)=P(e^{ix}),$$ hence $Q=P$.  Malik's inequality therefore gives:
$$2|P'(e^{ix})|\leq 2n, \hbox{\quad that is}\quad |P'(e^{ix})|=\sqrt{(T'(x))^2+n^2 (T(x))^2}\leq n.$$
\end{demo}

\noindent An obvious corollary is once more

\begin{theoreme}[Bernstein.] Let  $T$ be a \textnormal{complex}  trigonometric polynomial of degree  $n$, with $|T(x)|\leq 1$ for all $x\in \R$.  Then
$$|T'(x)|\leq n \hbox{\quad for all }\quad  x\in \R. $$

\end{theoreme}
\bigskip

\begin{demo}let $u$ be a unimodular  complex number  and $S_u=\Re (uT)$, a real trigonometric  polynomial of degree $n$ satisfying $|S_{u}(x)|\leq 1$ for all $x\in \R$. By the Schaake-van der Corput theorem,  we have
$| \Re (u T'(e^{ix}))|\leq n$, whence the result, optimizing with respect to  $u$. \end{demo}
\bigskip

\section{Bernstein's inequality via integral representation}

In this section we provide an approach to Bernstein's inequality for
the sup-norm, the $L^{p}-$norm ($p\geq1$) and some other variants,
based on new integral representations for algebraic/trigonometric
polynomials. The latter are developed in \cite{BaZa1,BaZa2} in a
more general context, to prove Bernstein-type inequalities for rational
functions. These integral representations are footed on the theory
of \textit{model spaces} and their reproducing kernels. The model
spaces are the subspaces of the Hardy space $H^{2}$ which are invariant
with respect to the backward shift operator, (we refer to \cite{Nik2}
for the general theory of model spaces and their numerous applications).
Applying this method to the case of algebraic polynomials, Bernstein's
inequalities for the sup-norm and for the $L^{p}-$norm ($p\geq1$)
are easily demonstrated. \\
 Our integral representations require to introduce the scalar product
$\left\langle \cdot,\,\cdot\right\rangle $ on $L^{2}=L^{2}(\mathbb{T})$

\[
\left\langle f,\,g\right\rangle =\int_{\mathbb{T}}f(u)\overline{g(u)}{\rm d}m(u).
\]
For $n\geq1$ the (algebraic) Dirichlet kernel $D_{n}$ is defined
as
\[
D_{n}(z)=\sum_{k=0}^{n-1}z^{k}.
\]

\subsection{\label{subsec:The-case-of-algebraic-poly}The case of algebraic polynomials}

Given an algebraic polynomial of degree $n$, $P(z)=\sum_{k=0}^{n}a_{k}z^{k}$
and given $\xi$ in the closed unit disk, we have
\[
P'(\xi)=\sum_{k=1}^{n}ka_{k}\xi^{k-1}=\bigg<P(z),\,z\frac{1}{(1-\overline{\xi}z)^{2}}\bigg>.
\]
Expanding $(1-(\overline{\xi}z)^{n})^{2}$ we observe that
\[
z\frac{1}{(1-\overline{\xi}z)^{2}}-z\frac{(1-(\overline{\xi}z)^{n})^{2}}{(1-\overline{\xi}z)^{2}}=z\frac{1}{(1-\overline{\xi}z)^{2}}-z\left(D_{n}(\overline{\xi}z)\right)^{2}
\]
is orthogonal to $P$. This yields
\begin{equation}
P'(\xi)=\int_{\mathbb{T}}P(u)\overline{u\left(D_{n}(\overline{\xi}u)\right)^{2}}{\rm d}m(u),\qquad|\xi|\leq1.\label{intrepder_anal_poly}
\end{equation}
Therefore for any unimodular $\xi$
\[
\left|P'(\xi)\right|\leq\left\Vert P\right\Vert _{\infty}\left\Vert D_{n}\right\Vert _{2}^{2}
\]
and Bernstein's inequality for the sup-norm follows. Following the
same approach we prove that
\[
\left\Vert P'\right\Vert _{p}\leq n\left\Vert P\right\Vert _{p},\qquad p\in[1,\infty].
\]
\begin{demo} An application of \eqref{intrepder_anal_poly} indeed
yields
\begin{eqnarray*}
\left\Vert P'\right\Vert _{p}^{p} & = & \int_{\mathbb{T}}\left|P'(\xi)\right|^{p}{\rm d}{m}(\xi)\\
 & = & \int_{\mathbb{T}}\bigg|\int_{\mathbb{T}}P(u)\overline{u\left(D_{n}(\overline{\xi}u)\right)^{2}}{\rm d}{m}(u)\bigg|^{p}{\rm d}{m}(\xi)\\
 & \leq & \int_{\mathbb{T}}\bigg(\int_{\mathbb{T}}\left|P(u)\right||D_{n}(\overline{\xi}u)|^{2}\bigg){\rm d}{m}(u)\bigg)^{p}{\rm d}{m}(\xi).
\end{eqnarray*}
We apply Hölder's inequality ($q$ is the conjugate exponent of $p:$
$\frac{1}{p}+\frac{1}{q}=1$)

\begin{multline*}
\bigg(\int_{\mathbb{T}}\left|P(u)\right||D_{n}(\overline{\xi}u)|^{2}{\rm d}{m}(u)\bigg)^{p}\\
\le\bigg(\int_{\mathbb{T}}|D_{n}(\overline{\xi}u)|^{2}{\rm d}{m}(u)\bigg)^{\frac{p}{q}}\int_{\mathbb{T}}\left|P(u)\right|^{p}|D_{n}(\overline{\xi}u)|^{2}{\rm d}{m}(u)\\
\leq n^{\frac{p}{q}}\int_{\mathbb{T}}\left|P(u)\right|^{p}|D_{n}(\overline{\xi}u)|^{2}{\rm d}{m}(u).
\end{multline*}
It remains to integrate with respect to $\xi$ and apply the Fubini-Tonelli
theorem to conclude. \end{demo}

The case of trigonometric polynomials is more technical and removed
to the end of the section. More precisely, in subsection \ref{subsec:The-case-of-trig-poly}
we provide an analog of \eqref{intrepder_anal_poly} for trigonometric
polynomials $T$ of degree at most $n$. Applying ``roughly'' the
above approach to $T$ yields $\left\Vert T'\right\Vert _{p}\leq2n\left\Vert T\right\Vert _{p}$
instead of $\left\Vert T'\right\Vert _{p}\leq n\left\Vert T\right\Vert _{p}$.
\\

In the next subsection we show that the continuous embeddings of some
Besov/Wiener algebras of analytic functions on $\mathbb{D}$, into
the algebra of bounded analytic functions, are invertible over the
set of algebraic polynomials of degree at most $n.$ We discuss the
asymptotic behavior of the respective embedding constants as $n\rightarrow\infty.$

\subsection{\label{subsec:The-case-of-analytic-poly}Inequalities for algebraic
polynomials in Besov/Wiener algebras}

We denote by $H^{\infty}$ the algebra of bounded analytic functions
on $\mathbb{D}$ i.e. the space of holomorphic functions $f$ on $\mathbb{D}$
such that $\left\Vert f\right\Vert _{\infty}<\infty.$ Given a Banach
algebra $X$ continuously embedded into $H^{\infty}$, we are interested
in inequalities of the type
\[
\left\Vert P\right\Vert _{X}\leq C_{X}(n)\left\Vert P\right\Vert _{\infty}
\]
holding for any algebraic polynomial $P$ of degree at most $n$.
The selected algebras $X$ below, are of particular interest for applications
in matrix analysis and operator theory, see \cite{Nik1} for more
details.
\begin{enumerate}
\item $B_{1,1}^{1}$ is the Besov algebra of analytic functions $f$ on
$\mathbb{D}$ such that
\[
\left\Vert f\right\Vert _{B_{1,1}^{1}}^{*}:=\int_{\mathbb{D}}\left|f''(u)\right|{\rm d}A(u)<\infty
\]
where ${\rm d}A$ stands for normalized Lebesgue measure on $\mathbb{D}$ and $\left\Vert \cdot \right\Vert _{B_{1,1}^{1}}^{*}$ is a semi-norm on $B_{1,1}^{1}$.
Vitse's functional calculus \cite{VP} shows that given a Banach Kreiss
operator $A,$ i.e. an operator $A$ satisfying the resolvent estimate
\[
\left\Vert (\lambda-A)^{-1}\right\Vert \leq C(\left|\lambda\right|-1)^{-1},\qquad\left|\lambda\right|>1,
\]
we have
\[
\left\Vert P(A)\right\Vert \leq2C\left\Vert P\right\Vert _{B_{1,1}^{1}}^{*}
\]
for every algebraic polynomial $P$.
\item $W$ is the analytic Wiener algebra of absolutely converging Fourier/Taylor
series, i.e. the space of all $f=\sum_{k\geq0}a_{k}z^{k}$ such that:
\[
\left\Vert f\right\Vert _{W}:=\sum_{k\geq0}\left|a_{k}\right|<\infty.
\]
It is easily verified that for any operator $A$ acting on a Banach
space, satisfying $\left\Vert A\right\Vert \leq1$, we have
\[
\left\Vert P(A)\right\Vert \leq\left\Vert P\right\Vert _{W}
\]
for every algebraic polynomial $P$.
\item $B_{\infty,1}^{0}$ is the Besov algebra of analytic functions $f$
in $\mathbb{D}$ such that
\[
\left\Vert f\right\Vert _{B_{\infty,1}^{0}}^{*}:=\int_{0}^{1}\left\Vert f'_{r}\right\Vert _{\infty}{\rm d}r<\infty
\]
where $f_{r}(z)=f(rz)$ and $\left\Vert \cdot \right\Vert _{B_{\infty,1}^{0}}^{*}$ is a semi-norm on $B_{\infty,1}^{0}$. Let $A$ be power bounded operator on a Hilbert
space: $\sup_{k\geq0}\left\Vert A^{k}\right\Vert =a<\infty$. Peller's
functional calculus \cite{Pel2} shows that
\[
\left\Vert P(A)\right\Vert \leq k_{G}a^{2}\left\Vert P\right\Vert _{B_{\infty,1}^{0}}^{*}
\]
for every algebraic polynomial $P$, where $k_{G}$ is an absolute
(Grothendieck) constant.
\end{enumerate}
Observe that the following continuous embeddings hold
\[
B_{1,1}^{1}\subset W\subset B_{\infty,1}^{0}\subset H^{\infty}
\]
see \cite{BeLo,Pee} or \cite[Sect. B.8.7]{Nik2}. It turns out that
the continuous embeddings $W\subset H^{\infty}$, $B_{\infty,1}^{0}\subset H^{\infty}$
and $B_{1,1}^{1}\subset H^{\infty}$ are invertible on the space of
complex algebraic polynomials of degree at most $n\geq1$. More precisely
we prove the following inequalities.

\begin{proposition} \label{Prop1}For any algebraic polynomial $P$
of degree at most $n$ the following inequalities hold

\begin{equation}
\left\Vert P\right\Vert _{W}\leq\sqrt{n+1}\left\Vert P\right\Vert _{\infty},\label{eq:trivial_inequ}
\end{equation}
the bound $\sqrt{n+1}$ being the best possible asymptotically as
$n\rightarrow\infty$, and
\begin{equation}
\left\Vert P\right\Vert _{B_{\infty,1}^{0}}^{*}\leq\left(\sum_{k=1}^{n-1}\frac{1}{2k+1}\right)\left\Vert P\right\Vert _{\infty}.\label{eq:logn_inequ}
\end{equation}
\end{proposition}

The asymptotic sharpness of $\ln n$ over the space of algebraic polynomials
of degree $n$ as $n\rightarrow\infty$, is an open question. Let
us recall a result by V. Peller \cite[Corollary 3.9]{Pel2}.

\begin{proposition}[Peller] \label{cor3.9}Let $A$ be a power
bounded operator on a Hilbert space. Then there exists a postive $M$
such that for any algebraic polynomial $P$ of degree $n$ the following
inequality holds
\[
\left\Vert P(A)\right\Vert \leq M\ln(n+2)\left\Vert P\right\Vert _{\infty}.
\]

\end{proposition} Indeed combining Peller's functional calculus \cite{Pel2}
with \eqref{eq:logn_inequ} we find
\[
\left\Vert P(A)\right\Vert \leq k_{G}a^{2}\left\Vert P\right\Vert _{B_{\infty,1}^{0}}^{*}\leq k_{G}a^{2}\left(\ln n+\gamma+o(1)\right)\left\Vert P\right\Vert _{\infty},
\]
where $a=\sup_{k\geq0}\left\Vert A^{k}\right\Vert ,$ $k_{G}$ is
an absolute (Grothendieck) constant, and $\gamma$ is the Euler constant.
The asymptotic sharpness of $\ln n$ in Proposition \ref{cor3.9}
is also an open question.

\noindent \begin{demo} [Proof of Proposition \ref{Prop1}] We
first prove \eqref{eq:trivial_inequ}. Given $P=\sum_{k=0}^{n}a_{k}z^{k}$,
Cauchy-Schwarz inequality yields
\[
\left\Vert P\right\Vert _{W}\leq\sqrt{n+1}\left\Vert P\right\Vert _{2}\leq\sqrt{n+1}\left\Vert P\right\Vert _{\infty}.
\]
Moreover, the bound $\sqrt{n+1}$ is asymptotically sharp as shown
for example by Kahane (\cite{KA})  at the beginning of his construction of ultraflat
polynomials, when he produces polynomials $P(z)=\sum_{k=0}^{n}a_{k}z^{k}$
with $|a_{k}|=1$ for $k=0,\dots,  n$ and $\Vert P\Vert_{\infty}\geq(1-\delta_{n})\sqrt{n+1}$
where $\delta_{n}\to0^{+}$.

\noindent Now we prove \eqref{eq:logn_inequ}. Applying \eqref{intrepder_anal_poly}
with $\zeta=rv$ and $v\in\mathbb{T}$ we find
\[
P'(rv)=\int_{\mathbb{T}}P(u)\overline{u\left(D_{n}(r\overline{v}u)\right)^{2}}{\rm d}m(u).
\]
This yields
\begin{align*}
\left|P'(rv)\right| & \leq\left\Vert P\right\Vert _{\infty}\int_{\mathbb{T}}\left|\overline{u\left(D_{n}(r\overline{v}u)\right)^{2}}\right|{\rm d}m(u)\\
 & =\left\Vert P\right\Vert _{\infty}\sum_{k=0}^{n-1}r^{2k}.
\end{align*}
Therefore taking the supremum over unimodular $\xi$ and integrating
over $r\in[0,1]$ we get
\[
\left\Vert P\right\Vert _{B_{\infty,1}^{0}}^{*}\leq\left\Vert P\right\Vert _{\infty}\sum_{k=0}^{n-1}\frac{1}{2k+1}.
\]

\end{demo}

We finally treat the case of the $B_{1,1}^{1}-$norm of $P$. Since
the second derivative of $P$ is involved in the definition of $\left\Vert P\right\Vert _{B_{1,1}^{1}}^{*}$
we first need to give an analog of \eqref{intrepder_anal_poly} for
$P''$. Clearly,
\[
P''(\xi)=\sum_{k=2}^{n}k(k-1)a_{k}\xi^{k-2}=2\bigg<P,\,z^{2}\frac{1}{(1-\overline{\xi}z)^{3}}\bigg>.
\]
Expanding $(1-(\overline{\xi}z)^{n})^{3}$ we observe that
\[
z^{2}\frac{1}{(1-\overline{\xi}z)^{3}}-z^{2}\frac{(1-(\overline{\xi}z)^{n})^{3}}{(1-\overline{\xi}z)^{3}}
\]
is orthogonal to any polynomial of degree at most $n+1$ and especially
to $P$. Therefore
\begin{equation}
P''(\xi)=2\int_{\mathbb{T}}P(u)\overline{u^{2}\left(D_{n}(\overline{\xi}u)\right)^{3}}{\rm d}m(u),\qquad|\xi|\leq1.\label{int_rep_2nd_der}
\end{equation}
We will use \eqref{int_rep_2nd_der} to prove next proposition.

\begin{proposition}[Vitse, Peller, Bonsall-Walsh] \label{Prop2}
For any algebraic polynomial $P$ of degree at most $n$ the following
inequality holds
\[
\left\Vert P\right\Vert _{B_{1,1}^{1}}^{*}\leq\frac{8}{\pi}\left(\sum_{k=0}^{n-1}\frac{\Gamma\left(k+\frac{3}{2}\right)^{2}}{k!(k+1)!}\right)\left\Vert P\right\Vert _{\infty},
\]
where $\Gamma$ is the standard Euler Gamma function. In particular
\begin{equation}
\left\Vert P\right\Vert _{B_{1,1}^{1}}^{*}<\frac{8}{\pi}n\left\Vert P\right\Vert _{\infty}.\label{eq:BW_inequ}
\end{equation}
\end{proposition} It is shown by P. Vitse in \cite[Lemma 2.3]{VP}
that \eqref{eq:BW_inequ} actually holds for rational functions $r$
of degree $n$ whose poles lie outside the closed unit disk, with
the same numerical constant $\frac{8}{\pi}$. Note that the same inequality
was originally proved by V. Peller in \cite{Pel1} without giving
an explicit numerical constant. Vitse's proof makes use of a theorem
by F. F. Bonsall and D. Walsh \cite{BoWa}, where the constant $\frac{8}{\pi}$
is sharp. The proof below does not make use of the theory of Hankel
operators, and is only based on \eqref{int_rep_2nd_der}.

\begin{demo}[Proof of Proposition \ref{Prop2}] We rewrite \eqref{int_rep_2nd_der}
as
\begin{equation}
P''(\xi)=2\left\langle P,\,z^{2}(1-(\overline{\xi}z)^{n})\left(\frac{1-(\overline{\xi}z)^{n}}{(1-\overline{\xi}z)^{\frac{3}{2}}}\right)^{2}\right\rangle ,\label{intrep}
\end{equation}
and  we use the Taylor expansion of $\left(1-\overline{\xi}z\right)^{-\frac{3}{2}}$ to get
\begin{eqnarray*}
\frac{1-(\overline{\xi}z)^{n}}{\left(1-\overline{\xi}z\right)^{\frac{3}{2}}} & = & (1-(\overline{\xi}z)^{n})\sum_{k\geq0}\frac{\Gamma\left(k+\frac{3}{2}\right)}{k!\Gamma\left(\frac{3}{2}\right)}\overline{\xi}^{k}z^{k}\\
 & = & \frac{2}{\sqrt{\pi}}(1-(\overline{\xi}z)^{n})\sum_{k\geq0}\frac{\Gamma\left(k+\frac{3}{2}\right)}{k!}\overline{\xi}^{k}z^{k}\\
 & = & \frac{2}{\sqrt{\pi}}\Big(\sum_{k\geq0}\frac{\Gamma\left(k+\frac{3}{2}\right)}{k!}\overline{\xi}^{k}z^{k}-\sum_{k\geq0}\frac{\Gamma\left(k+\frac{3}{2}\right)}{k!}\overline{\xi}^{k+n}z^{k+n}\Big)\\
 & = & \frac{2}{\sqrt{\pi}}(\varphi_{\xi}(z)-\psi_{\xi}(z)),
\end{eqnarray*}
where $\varphi_{\xi}(z)=\sum_{k\geq0}\frac{\Gamma\left(k+\frac{3}{2}\right)}{k!}\overline{\xi}^{k}z^{k}$
and $\psi_{\xi}(z)=\sum_{k\geq0}\frac{\Gamma\left(k+\frac{3}{2}\right)}{k!}\overline{\xi}^{k+n}z^{k+n}.$
We observe that the functions
\[
z\mapsto z^{2}(1-(\overline{\xi}z)^{n})(\psi_{\xi}(z))^{2}
\]
and
\[
z\mapsto z^{2}(1-(\overline{\xi}z)^{n})\varphi_{\xi}(z)\psi_{\xi}(z),
\]
are orthogonal to any algebraic polynomial of degree at most $n$.
Writing
\[
\left(\frac{1-(\overline{\xi}z)^{n}}{\left(1-\overline{\xi}z\right)^{\frac{3}{2}}}\right)^{2}=\frac{4}{\pi}((\varphi_{\xi}(z))^{2}+(\psi_{\xi}(z))^{2}-2\varphi_{\xi}(z)\psi_{\xi}(z)),
\]
\eqref{intrep} becomes
\begin{eqnarray*}
P''(\xi) & = & 2\left\langle P,\,z^{2}(1-(\overline{\xi}z)^{n})\left(\frac{1-(\overline{\xi}z)^{n}}{(1-\overline{\xi}z)^{\frac{3}{2}}}\right)^{2}\right\rangle \\
 & = & \frac{8}{\pi}\left\langle P,\,z^{2}(1-(\overline{\xi}z)^{n})((\varphi_{\xi}(z))^{2}+(\psi_{\xi}(z))^{2}-2\varphi_{\xi}(z)\psi_{\xi}(z))^{2}\right\rangle \\
 & = & \frac{8}{\pi}\left\langle P,\,z^{2}(1-(\overline{\xi}z)^{n})(\varphi_{\xi}(z))^{2}\right\rangle \\
 & = & \frac{8}{\pi}\left\langle P,\,(z\varphi_{\xi}(z))^{2}\right\rangle .
\end{eqnarray*}
since $z\mapsto z^{n+2}(\varphi_{\xi}(z))^{2}$ is also orthogonal
to $P$. Denoting by $S_{n}^{\xi}(z)=\sum_{k=0}^{n-1}\frac{\Gamma\left(k+\frac{3}{2}\right)}{k!}\overline{\xi}^{k}z^{k}$
and $R_{n}^{\xi}(z)=\sum_{k\geq n}\frac{\Gamma\left(k+\frac{3}{2}\right)}{k!}\overline{\xi}^{k}z^{k},$we
have
\[
S_{n}^{\xi}(z)+R_{n}^{\xi}(z)=\varphi_{\xi}(z),
\]
and
\begin{eqnarray*}
(z\varphi_{\xi}(z))^{2} & = & (zS_{n}^{\xi}(z)+zR_{n}^{\xi}(z))^{2}\\
 & = & (zS_{n}^{\xi}(z))^{2}+(zR_{n}^{\xi}(z))^{2}+2z^{2}R_{n}^{\xi}(z)S_{n}^{\xi}(z),
\end{eqnarray*}
where the two last terms are again orthogonal to $P.$ Finally, we
obtain the following integral representation:
\[
P''(\xi)=\frac{8}{\pi}\left\langle P,\,(zS_{n}^{\xi}(z))^{2}\right\rangle .
\]
Using the standard Cauchy duality we have that for any $\xi\in\mathbb{D}$,
\[
\left|P''(\xi)\right|\leq\frac{8}{\pi}\left\Vert P\right\Vert _{\infty}\int_{\mathbb{T}}\left|uS_{n}^{\xi}(u)\right|^{2}{\rm d}m(u).
\]
Integrating over $\mathbb{D}$ with respect to the normalized area
measure, we find
\begin{align*}
\int_{\mathbb{D}}\left|P''(\xi)\right|{\rm d}A(\xi) & =\frac{8}{\pi}\left\Vert P\right\Vert _{\infty}\int_{\mathbb{D}}\left|\int_{\mathbb{T}}\overline{(uS_{n}^{\xi}(u))^{2}}{\rm d}m(u)\right|{\rm d}A(\xi)\\
 & \leq\frac{8}{\pi}\left\Vert P\right\Vert _{\infty}\int_{\mathbb{D}}\left(\int_{\mathbb{T}}\left|S_{n}^{\xi}(u)\right|^{2}{\rm d}m(u)\right){\rm d}A(\xi)\\
 & =\frac{8}{\pi}\left\Vert P\right\Vert _{\infty}\int_{\mathbb{T}}\left(\int_{\mathbb{D}}\left|S_{n}^{\xi}(u)\right|^{2}{\rm d}A(\xi)\right){\rm d}m(u).
\end{align*}
We conclude noticing that $\int_{\mathbb{D}}\left|S_{n}^{\xi}(u)\right|^{2}{\rm d}A(\xi)$
is the square of norm of $S_{n}^{\xi}$ in the standard Bergman space,
we find
\[
\int_{\mathbb{D}}\left|S_{n}^{\xi}(u)\right|^{2}{\rm d}A(\xi)=\sum_{k=0}^{n-1}\frac{\Gamma\left(k+\frac{3}{2}\right)^{2}}{(k+1)(k!)^{2}}\left|u\right|^{2k}=\sum_{k=0}^{n-1}\frac{\Gamma\left(k+\frac{3}{2}\right)^{2}}{k!(k+1)!}.
\]
\end{demo}

\subsection{\label{subsec:The-case-of-trig-poly}The case of trigonometric polynomials}

As a generalization of \eqref{intrepder_anal_poly} we prove the following
integral representation for the derivative of trigonometric polynomials.

\begin{lemme}\label{Lemma_int_rep_trig_poly} For any trigonometric
polynomial $T$ of degree $n$ we have
\begin{equation}
T'(\xi)=\left\langle T,\,K_{\xi}\right\rangle ,\qquad|\xi|=1,\label{integrepder}
\end{equation}
where for all $u,\,\xi\in\mathbb{T},$ $K_{\xi}(u)=u\left(D_{n}(\overline{\xi}u)\right)^{2}-\xi^{2}\overline{u}\overline{\left(D_{n}(\overline{\xi}u)\right)^{2}}.$

\end{lemme} The proof of Bernstein's inequality for $p\in[1,\infty]$
with constant $2n$ instead of $n$, follows from the above lemma.
Indeed, following the same trick as in subsection \ref{subsec:The-case-of-algebraic-poly}
we get
\[
\left\Vert T'\right\Vert _{p}\leq2n\left\Vert T\right\Vert _{p},\qquad p\in[1,\infty].
\]

\begin{demo}{[}Proof of the lemma of integral representation{]} We
put $T=\sum_{k=-n}^{n}a_{k}z^{k}$ ($z=e^{it}$), $P=\sum_{k=0}^{n}a_{k}z^{k}$
and $R=\sum_{k=-m}^{-1}a_{k}z^{k}$ so that $T=P+R.$ Applying \eqref{intrepder_anal_poly}
to the algebraic polynomial $P$ we get
\begin{equation}
P'(\xi)=\bigg<P(z),\,z\left(D_{n}(\overline{\xi}z)\right)^{2}\bigg>=\bigg<T(z),\,z\left(D_{n}(\overline{\xi}z)\right)^{2}\bigg>,\label{eq:}
\end{equation}
because $R\perp z\left(D_{n}(\overline{\xi}z)\right)^{2}.$ We will
now perform a similar task for $R.$ Consider for this the algebraic polynomial
\[
Q(z)=\bar{z}R\left(\bar{z}\right),\qquad \bar{z}=1/z,
\]
whose degree does not exceed $n-1$. For the reasons given above,
we have
\[
Q(\xi)=\bar{\xi}R\left(\bar{\xi}\right)=\bigg<Q(z),\,\frac{1}{1-\overline{\xi}z}\bigg>,
\]
that is to say
\[
R(\frac{1}{\xi})=\xi\bigg<Q(z),\,\frac{1}{1-\overline{\xi}z}\bigg>.
\]
Deriving again with respect to $\xi$ we get
\[
-\frac{1}{\xi^{2}}R'(\frac{1}{\xi})=\bigg<Q(z),\,\frac{1}{1-\overline{\xi}z}\bigg>+\xi\bigg<Q(z),\,\frac{z}{(1-\overline{\xi}z)^{2}}\bigg>,
\]
that is to say
\[
R'(\frac{1}{\xi})=\bigg<Q(z),\,\frac{-\overline{\xi}^{2}}{(1-\overline{\xi}z)^{2}}\bigg>=\bigg<Q(z),\,\frac{-\overline{\xi}^{2}(1-\overline{\xi}^{n}z^{n})^{2}}{(1-\overline{\xi}z)^{2}}\bigg>,
\]
the last equality being due to the fact that $\frac{1}{(1-\overline{\xi}z)^{2}}-\frac{(1-(\overline{\xi}z)^{n})^{2}}{(1-\overline{\xi}z)^{2}}$
is orthogonal to any algebraic polynomial of degree at most $n$.
Rewriting this last equality using the integral representation of
the scalar product, we find
\[
R'(\frac{1}{\xi})=-\int_{\mathbb{T}}\bar{u}R\left(\bar{u}\right)\frac{\xi^{2}(1-\xi^{n}\bar{u}^{n})^{2}}{(1-\xi\bar{u})^{2}}{\rm d}m(u).
\]
Performing the variable change $v=\bar{u}$ and replacing $\xi$ with
$\bar{\xi}$ we obtain
\[
R'(\xi)=-\int_{\mathbb{T}}R\left(v\right)\frac{v\xi^{2}(1-\bar{\xi}^{n}v^{n})^{2}}{(1-\bar{\xi}v)^{2}}{\rm d}m(v).
\]
Finally
\begin{equation}
R'(\xi)=\bigg<R(z),\,-\xi^{2}\overline{z}\overline{\left(D_{n}(\overline{\xi}z)\right)^{2}}\bigg>=\bigg<T(z),\,-\xi^{2}\overline{z}\overline{\left(D_{n}(\overline{\xi}z)\right)^{2}}\bigg>,\label{eq:-1}
\end{equation}
because $P\perp\overline{z}\overline{\left(D_{n}(\overline{\xi}z)\right)^{2}}$.
It remains to combine \eqref{eq:} and \eqref{eq:-1} to complete
the proof.

\end{demo}
\medskip

As we will now see, subharmonicity and complex methods allow us to go beyond convexity and to consider $L^p$-quasi-norms for $0<p<1$, even for $p=0$ (the Mahler norm). Indeed, we begin with the Mahler norm.
\section{Case $p=0$, Mahler's result}
\noindent This section and the next one owe much to conversations with F.~Nazarov (\cite{NAZ}). Before Nazarov, the possible use of subharmonicity was alluded to  by the referee of Mahler's paper. But to our knowledge, none of the approaches that  follow theorem \ref{mah}  was detailed anywhere in the literature.\\
We will first show, following Mahler (\cite{MAH2}), that
\begin{theoreme}[Mahler.]\label{mah} It holds
 \begin{equation}\label{mahl}\Vert P'\Vert_0\leq n\Vert P\Vert_0\end{equation}
  for every \textnormal{algebraic}  polynomial  $P(z)=\sum_{k=0}^n a_k z^k$ of degree $n$.
  \end{theoreme}
\begin{demo}the proof, simpler than Mahler's initial one, consists of two steps.\\

 \noindent {\bf 1.} The result holds true if all roots of $P$ lie in $\overline{\D}$. Indeed, the same holds for $P'$ (by Gauss-Lucas) and in that case both members of the inequality (\ref{mahl}) are equal to $n|a_n|$, by Jensen's formula.\\

\noindent{\bf 2.} The result holds true in the general case. To see that, write
$$P(z)=a_n \prod_{|z_j|<1}(z-z_j)\,\prod_{|z_j|\geq1}(z-z_j)$$
$$Q(z)=a_n \prod_{|z_j|<1}(z-z_j)\,\prod_{|z_j|\geq1}(1-\overline{z_j}z).$$
All roots of $Q$ lie in $\overline{\D}$, and  $|P(z)|=|Q(z)|$ for $|z|=1$, so $|P'(z)|\leq|Q'(z)|$ for $|z|=1$ by Lemma   \ref{duo}. The first step now implies
$$\Vert P'\Vert_0\leq\Vert Q'\Vert_0\leq n\Vert Q\Vert_0= n \Vert P\Vert_0.$$
It is convenient to ``stock'' the result under the form:
\begin{equation}\label{eins} \int \log |P'(z)/n|dm(z)\leq  \int \log |P(z)|dm(z).\end{equation}
\end{demo}
\medskip
We will now show that, more generally (it seems that Mahler only treated the algebraic case):

\begin{theoreme}\label{mg} One has
 \begin{equation}\label{mn}\Vert T'\Vert_0\leq n\Vert T\Vert_0\end{equation}
  for every \textnormal{trigonometric}  polynomial  $T(z)=\sum_{k=-n}^n a_k z^k$ of degree $n$.
    \end{theoreme}
  \begin{demo}we first observe the following: if we write $T(x)=\sum_{k=-n}^n a_k e^{ikx}$, then $T'(x)=izT'(z)$ when $z=e^{ix}$ and $|T'(x)|=|T'(z)|$. We can thus work indifferently with the variable $x$ or the variable $z$ to prove our inequality.  Denote  $Q(z)=z^n T(z)$,  an algebraic polynomial of degree $2n$. Write

$$Q(z)=c\prod_{j=1}^{2n} (z-z_j)$$ where  $z_1,\ldots,z_p$ denote the roots of  modulus $\leq 1$, and  $z_{p+1},\ldots, z_{2n}$ those of  modulus $>1$ if some  exist. One has:

$$z\frac{T'(z)}{T(z)}=z\frac{Q'(z)}{Q(z)}-n=\sum_{j=1}^{2n}\frac{z}{z-z_j}-n $$
so that
\begin{equation}\label{etap1} \int_{\T}\log |T'(z)/T(z)|dm(z)=\int_{\T} \log \Big|\sum_{j=1}^{2n}\frac{z}{z-z_j}-n\Big|dm(z)=:M(z_1,\ldots,z_p)  \end{equation}
where $M$ is the function of   $p$ complex variables defined by
$$M(Z_1,\ldots, Z_p)=\int_{\T} \log \Big|\sum_{j=1}^{p}\frac{z}{z-Z_j}+\sum_{j=p+1}^{2n}\frac{z}{z-z_j}-n\Big|dm(z).$$
To emphasize the key role of subharmonicity, we first outline the following
\begin{lemme}\label{subh}One considers the two  functions
$$M(w)=\log \Big|\frac{1}{w-u}+v\Big|,\quad N(w)=\int_{\T}\log \Big|\frac{1}{w-z}+h(z)\Big|dm(z)$$
where $u\in \T$ and $v\in \C$, and where $h$ is a continuous  function on $\T$. Then
$M=:M_{u,v}$ is sub-harmonic on $\D$ and $N$ sub-harmonic on $\D$  and continuous on $\overline{\D}$.

\end{lemme}
\begin{demo}for $M=M_{u,v}$, it is enough to remark that it is the logarithm of  $|f|$ where  $f(w)=\frac{1}{w-u}+v$ is a holomorphic function on $\D$ since $|u|=1$. Next,
$$N=\int_{\T} M_{z, h(z)}dm(z)$$
(vector-valued integral) and a sum  (or a  barycenter) of subharmonic functions is again subharmonic. The  continuity of $N$ on $\overline{\D}$ results from classical integration theorems  (uniform integrability).
\end{demo}

\noindent The proof of Theorem \ref{mn} is then split into two steps:\\

\noindent {\bf 1.} One can assume  $|z_j|=1$ for $1\leq j\leq p$. Indeed, $M$ has the form:
$$M(Z_1,\ldots,Z_p)=\int_{\T} \log \Big|\sum_{j=1}^{p}\frac{z}{z-Z_j}+\varphi(z)\Big|dm(z)$$
where $\varphi$ is a fixed  continuous function on the  circle $\T$, hence by Lemma \ref{subh},  $M$ is separately sub-harmonic in $\D^p$ and separately  continuous in $\overline{\D}^p$.  Repeatedly applying to it the  maximum principle in one variable, one sees that there exist $(u_1,\ldots, u_p)\in \partial{\D}^p$,  the distinguished boundary of  $\D^p$, such that:
$$M(z_1,\ldots, z_p)\leq M(u_1,\ldots, u_p).$$ It is thus enough to prove that $M(u_1,\ldots, u_p)\leq \log n$, with
$u_1,\ldots,u_p$ in place of the roots $z_1,\ldots,z_p$ of $Q$. In the sequel, we shall henceforth assume, without loss of generality, that  those roots satisfy
$$1\leq j\leq p\Rightarrow |z_j|=1 \hbox{\ and}\ p+1\leq j\leq 2n\Rightarrow |z_j|>1.$$
In particular, all roots $z_j$ of $Q$ have modulus $\geq 1$.\\

\noindent {\bf 2.} One has the  implication (an essential remark)
$$ |z|<1\Rightarrow \Re \frac{z}{z-z_j}<\frac{1}{2} \hbox{\quad for all}\ 1\leq j\leq 2n.$$
Indeed, an easy computation gives
$$\Re \Big(\frac{1}{2}-\frac{z}{z-z_j}\Big)=\frac{1}{2}\frac{|z_j|^2-|z|^2}{|z-z_j|^2}\cdot$$
It follows that, setting $f(z)=z\frac{Q'(z)}{Q(z)}-n=\sum_{j=1}^{2n}\frac{z}{z-z_j}-n$, one has for $|z|<1$:
 $$\Re f(z)\leq \sum_{j=1}^{2n}\Re \frac{z}{z-z_j} -n< \sum_{j=1}^{2n} \frac{1}{2}-n=0.$$
 The holomorphic function $f$ hence has no zeros in $\D$ and we are in the cases of equality in  Jensen's formula~:
$$ \int_{\T} \log |f_r|dm=\log |f(0)|=\log n$$
where  $r<1$ and $f_{r}(z)=f(rz)$.The rational fraction $f$ being ``well-behaved'', we let $r$ tend to  $1$ to get
$$\int_{\T} \log |f|dm=\log n $$
which implies via (\ref{etap1})
 $$\int_{\T}\log |T'(z)/T(z)|dm(z)= \int_{\T} \log |f|dm=\log n$$
 or again $ \Vert T'\Vert_0= n\,\Vert T\Vert_0$, and we are even in the cases of  equality in  Bernstein-Arestov's inequality when all roots have modulus $\geq 1$.
 \end{demo}

\section{Case $0<p<1$,  Arestov's result}

We will prove (the case $p\geq 1$ being already treated in the section ``Convexity'', but being recovered here as well) the following theorem, due to Arestov (\cite{ARE}):
\begin{theoreme}[Arestov.] Let $p>0$. It holds, for any trigonometric polynomial $T$ of degree $n$:
$$\Vert T'\Vert_p\leq n\Vert T\Vert_p.$$
\end{theoreme}
\begin{demo}instead of starting from Bernstein's inequality for $L^p,\ p=\infty$, and of generalizing,  one starts from  Bernstein's inequality for $L^0$ (initially due to Mahler; cf. \cite{MAH} and also the remark of the referee of Mahler's paper on the   maximum principle for subharmonic functions of several  variables) and  one goes up. This is done in two steps, each of which  uses an integral representation, in the style of Section 4.\\
 {\bf 1.} It holds
\begin{equation}\label{facil}  \log^{+} |v|=  \int_{\T} \log |v+uw|\,dm(w) \quad   \forall u\in \T, \quad\forall v\in \C.\end{equation}
Indeed, one can assume  $u=1$, given the translation-invariance of   $m$; then, one separates the cases $|v|\geq 1,\  |v|<1$, and one is always back to
$$\int_{\T} \log|1+aw|\,dm(w)=0 \hbox{\ if }\ |a|<1$$
which is nothing but the  harmonicity of $w\mapsto \log|1+aw|$ in $\D$.\\
We first note that (\ref{facil}) implies:
\begin{equation}\label{easy} \int_{\T}  \log^{+}  |T'(z)/n|dm(z)\leq \int_{\T}  \log^{+}  |T(z)|dm(z).\end{equation}
Indeed, for fixed $w\in \T$, one applies  (\ref{eins}) to the polynomial $T+w E_n$ with $E_{n}(z)=z^n$. One gets, since  $E'_n=nE_{n-1}$:
$$\int \log |T'(z)/n+w E_{n-1}(z)|dm(z)\leq  \int_{\T} \log |T(z)+w E_{n}(z)|dm(z).$$
One next integrates both members with respect to $w$,  uses  Fubini's theorem and applies the identity  (\ref{facil}) for fixed $z$ with $u=E_{n-1}(z),\  v=T'(z)/n$,  or with $u=E_{n}(z),\ v=T(z)$, to obtain  (\ref{easy}).\\

\noindent   {\bf 2.} It holds for $p>0$ and $u\geq 0$:
  \begin{equation}\label{facile} u^p=  \int_{0}^\infty \log^{+}(u/a)\,p^{2} a^{p-1}da.\end{equation}
Indeed, write the right-hand side  as $I=\int_{0}^{u}\log(u/a)\,p^{2} a^{p-1}da$, and  integrate by parts, differentiating the log, to get
$I=\int_{0}^{u} pa^{p-1}da=u^p.$\\
 Write $d\mu(a)=p^2 a^{p-1}da$ to save notation ($\mu$ depends on $p$). Identity (\ref{facile}) and Fubini used twice give, using also (\ref{easy}):
$$\int_{\T} |T'(z)/n|^{p}dm(z)=\int_{\T} \Big[\int_{0}^\infty \log^{+}( |T'(z)|/na) d\mu(a)\Big]dm(z)$$$$=
\int_{0}^\infty \Big[\int_{\T} \log^{+}( |T'(z)|/na) dm(z)\Big]d\mu(a)$$$$\leq \int_{0}^\infty \Big[\int_{\T} \log^{+}( |T(z)|/a) dm(z)\Big]d\mu(a)$$$$= \int_{\T} \Big[\int_{0} ^\infty\log^{+}( |T(z)|/a) d\mu(a)\Big]dm(z)=\int_{\T} |T(z)|^p\,dm(z).$$
This ends the proof of Arestov's theorem. \end{demo}

\medskip

\section{ Final Remarks}

 \noindent  {\bf 1.} Passing to the limit, when $n\to \infty$, in the Riesz relation (\ref{marcelor}) $$\frac{1}{2n^2}\sum_{r=1}^{2n}\frac{1}{\sin^{2}\frac{(2r-1)\pi}{4n}}=1$$
 easily gives the Euler formulas
 $$\sum_{r=1}^{\infty}(2r-1)^{-2}=\pi^{2}/8 \hbox{\quad and}\   \sum_{r=1}^{\infty}r^{-2}=\pi^{2}/6.$$

\noindent {\bf 2.} The Bernstein-Arestov inequalities for the  trigonometric polynomials $T_{n}(x)=\sum_{k=-n}^n a_k e^{ikx}$, namely
$$\Vert T'_n\Vert_p\leq n\Vert T_n\Vert_p$$ thus hold for all  $p\geq 0$ (\cite{ARE}).   A striking aspect of those inequalities is that the full question was still open in  1980, even for algebraic polynomials, in spite of partial nice contributions due to  Mat\'e and Nevai (\cite{MANE}), which appeared in  Annals of Math.! The authors prove that, for $0<p<1$ and $P$ an algebraic polynomial of degree $n$, it holds
$$\Vert P'\Vert_p \leq n(4e)^{1/p}\Vert P\Vert_p.$$

\noindent {\bf 3.}  The result of   (\cite{ARE}) is more precise: if $\chi:\R^+\to \R^+$ is increasing, differentiable  with $x\,\chi'(x)$ increasing as well, for example if $\chi(x)=x^p$ with $p>0$ or $\chi(x)=\log x$, one has
\begin{equation}\label{mieux} \int_{\T} \chi(|T'_{n}(z)|)dm(z)\leq  \int_{\T} \chi(|n\,T_{n}(z)|)dm(z).\end{equation}

\noindent {\bf 4.}  In the  case $p=\infty$,  the quite interesting proof of M.~Riesz (\cite{Riesz}) could inspire for a proof of the existence  of the  function $\varphi$ in Lemma \ref{zg}. This Riesz formula gives a more precise result than Bernstein's one, as we saw: \begin{equation}\label{plus} \vert T'(0)\vert \leq n  \sup_{x\in E_{n}}|T(x)|.\end{equation}
where  $E_n$ (a coset) designates the fixed subset of cardinality $2n$ formed by the numbers  $\frac{(2r-1)\pi}{2n},\quad 1\leq r\leq 2n$. More precisely, identifying   $\frac{k\pi}{2n}$ and $e^{i\frac{k\pi}{2n}}$,
 let $G_{4n}$ be the group of  $4n$-th roots of unity  and $H_n=G_{4n}^2$ be the  subgroup of order  $2n$ formed by squares. One has  $E_n=\omega H_n$. \\

\noindent  {\bf 5.} In (\cite{QUSA}), one can find various improvements of Bernstein's inequality for the so-called ultraflat polynomials $P$ of Kahane (those of degree $n$ with unimodular coefficients and with modulus nearly $\sqrt n$ on the unit circle), under the form
 $$\Vert P'\Vert_p \leq \gamma_{p}\, n \big(1+O(n^{-1/7})\big)\Vert P\Vert_p$$
 where $\gamma_p<1$ is a constant given in explicit terms.
 \bigskip

\section{Acknowledgments}

The first named author warmly thanks F.~Nazarov for very useful exchanges about Bernstein's inequality. The second-named  author acknowledges the Russian Science Foundation grant 14-41-00010. We finally thank the referee for his careful reading of our manuscript and his many insightful comments and suggestions.

\end{document}